\documentclass[preprint,12pt]{elsarticle_arxiv}

\usepackage{amssymb}

\usepackage{caption}
\usepackage{subcaption}

\usepackage{tabularx}
\usepackage{graphicx}
\usepackage{adjustbox}
\usepackage{url}
\usepackage{float}

\journal{Renewable Energy}

\begin{document}

\begin{frontmatter}



\title{Sustainable Hydropower Planning in Gabon}


\author[inst1]{Rafael Kelman}
\author[inst1]{Tainá Martins Cunha}
\author[inst1]{Luiz Rodolpho Sauret Cavalcanti de Albuquerque}
\author[inst1]{Marcelo Gomes Metello}
\author[inst1]{Tarcisio Luiz Coelho de Castro}

\affiliation[inst1]{organization={PSR},
            addressline={Praia de Botafogo, 370 - Botafogo}, 
            city={Rio de Janeiro},
            postcode={22250-040}, 
            state={Rio de Janeiro},
            country={Brazil}}

\begin{abstract}
Hydropower is a renewable, controllable, and flexible source of electricity. These are instrumental features to support decarbonization efforts, as an enabler of non-controllable and variable sources of renewable electricity. Sometimes hydropower is accompanied by other services provided by multipurpose reservoirs, such as water supply, irrigation, navigation, flood control and recreation. Despite all these benefits, hydropower can be a polarizing issue. A large sample of projects with poor planning and execution provides numerous arguments for its opponents. Large and complex projects frequently suffer overcost and delays. The direct impacts are related to the disruption of river ecosystems and surrounding habitats due to the flooding of large areas, and the fragmentation of rivers caused by the construction of dams and the reduction of sediment transport that impoverishes aquatic life. People displacement and compensation are always a complex issue. Hydropower projects can also cause indirect impacts, such as additional deforestation related to the construction of workers’ villages, access roads and transmission lines. Finally, reservoirs may also become a significant source of methane emission, especially in tropical areas. This paper offers an analytical approach for sustainable hydropower planning at the river basin scale called HERA: Hydropower and Environmental Resource Assessment. Built on three main components, namely geoprocessing, engineering, and optimization, HERA screens and compares hydropower development alternatives to guarantee social and environmental objectives while maximizing mentioned economic benefits. It was designed to encourage a transparent and participatory hydropower planning process from the early stages. Experience has shown this approach increases the chances of better and balanced outcomes. A case study is presented in the Ogooue river basin in Gabon.
\end{abstract}

\begin{graphicalabstract}
\centering
\includegraphics[width=1\textwidth]{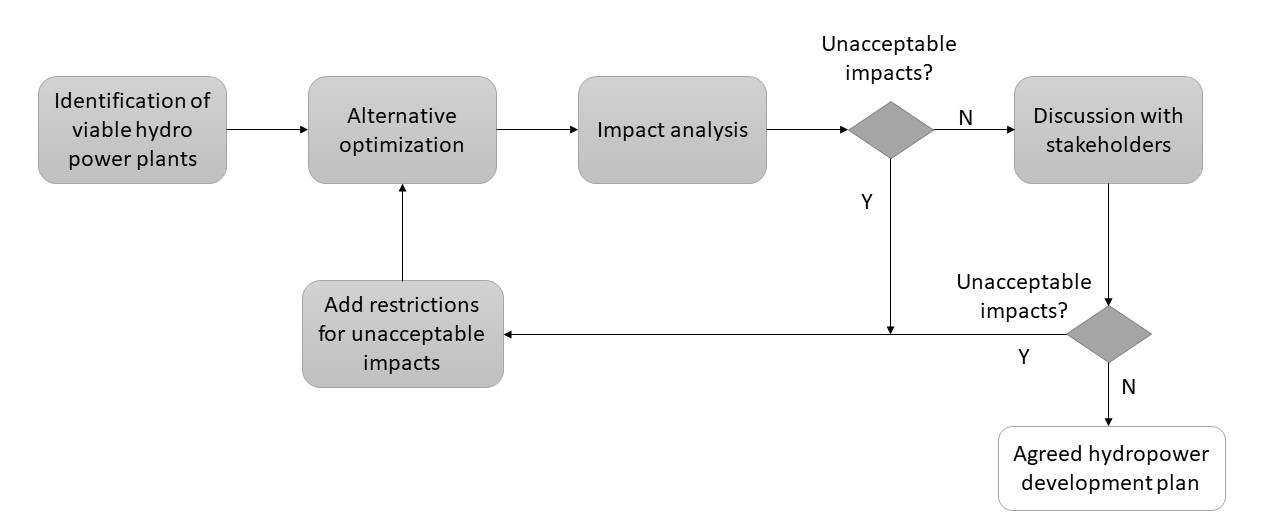}
\end{graphicalabstract}



\end{frontmatter}


\section{Introduction}
\label{1}
Hydropower is generally viewed as a proven infrastructure development path for developing countries. The need for reliable electric power is closely linked to economic growth \cite{1,2}. Moreover, in the current global context of energy transition, even developing countries are increasingly looking to make their development decisions that are in line with principles of sustainability \cite{3,4}. They may be able to leap-frog over older generation, classical hydropower dam development. Hydroelectric plants have unique characteristics that can provide important services to power grids in addition to electricity production, such as fast response times and operating reserves, that help control system frequency for changes in supply and demand. These services can support the development of variable renewable energy (VRE) sources globally, such as wind and solar power \cite{5,6}.

Poorly planned or badly executed hydroelectric projects, however, with frequent cost overruns and long delays, lack of transparency and insufficient communications with local communities, have increased the opposition to hydropower \cite{7,8,9}. Public concerns over large social and environmental impacts, including land flooding, migration barriers, loss of biodiversity, changes to flow regimes and so on \cite{10}. These impacts are especially  important considering that a significant share of the remaining hydropower potential lies in the world’s most biodiverse river basins, such as the Amazon, Congo, and Mekong \cite{11}. If it is clear the hydropower development will occur in a river basin, then the selection of which dams in which locations becomes highly relevant for the conservation of biodiversity \cite{11}, minimization of trade-offs between power generation, sediment supply and nutrient transport \cite{12}, and the minimization of impacts to local communities \cite{13}. Thus, well-considered socioenvironmental indicators must be determined from the early stages of the hydropower development planning process, that is, well before dam sites are chosen \cite{14}. This facilitates an objective trade-off analysis of electricity production, social impacts, and conservation to be made by stakeholders and government planners in a genuinely transparent decision-making process. The benefits of this process are an increase in the chance of a credible and successful negotiation process that minimizes risks for project developers and investors \cite{15}, as well as civil society stakeholders \cite{16}.

Several studies have demonstrated how site-specific assessment protocols can largely ignore cumulative impacts, while strategic, basin-wide planning enhances the probability of selecting dam configurations with less destructive, aggregate environmental outcomes \cite{17, 18, 19}. Proposed dam sites must be evaluated within the context of sustaining a portfolio of ecosystem services and biodiversity conservation, and alternative sites should be considered explicitly \cite{11,20,21} 

One of the major concerns about hydropower development in freshwater ecosystems is river fragmentation \cite{22,23}. Proposed dams can block the migratory patterns of fishes and threaten species with habitat loss and fragmentation \cite{24}. In time, such habitat alterations can reduce fish populations, impact local fisheries, and cause local extinctions \cite{4}. Dams also present the threat of sediment entrapment, resulting in alteration of river channel and floodplain geomorphology and associated ecosystem services. Despite important threats to biodiversity and sediment regimes, river fragmentation is seldom considered in hydropower planning in practice. 

New analytical tools and high-resolution environmental data can clarify trade-offs between engineering and environmental goals and can enable governments and funding institutions to compare alternative sites for dam building \cite{18,20}. In this paper we introduce a methodology that we refer to as Hydropower and Environmental Resource Assessment (HERA), accompanied by a computational model that applies this methodology. HERA simulates the construction of numerous projects in different candidate sites, quantifies their impacts and shortlists the alternative (set of projects) that maximizes the combined economic benefit considering socioenvironmental targets or constraints. 

Planning for sustainable hydropower development in Gabon illustrates the application of this methodology. Gabon is located in the Atlantic coast of Central Africa. It is one of the world’s few net sequesters of carbon \cite{25} with 85\% of its land covered in carbon-absorbing rainforest - an area about the size of the UK. The population of just 2 million is 90\% urban and the absence of major highways helps maintain the country's forests in pristine conditions. Its forests are part of the Congo Basin rainforest, the planet’s most important forest ecosystem after the Amazon \cite{25}. But Gabon’s biodiversity relevance goes far beyond forests. With over 350 species of fresh and brackish water fishes, the country also represents a freshwater biodiversity hotspot \cite{26}.

\begin{figure}[h]
    \centering
    \includegraphics[width=1\textwidth]{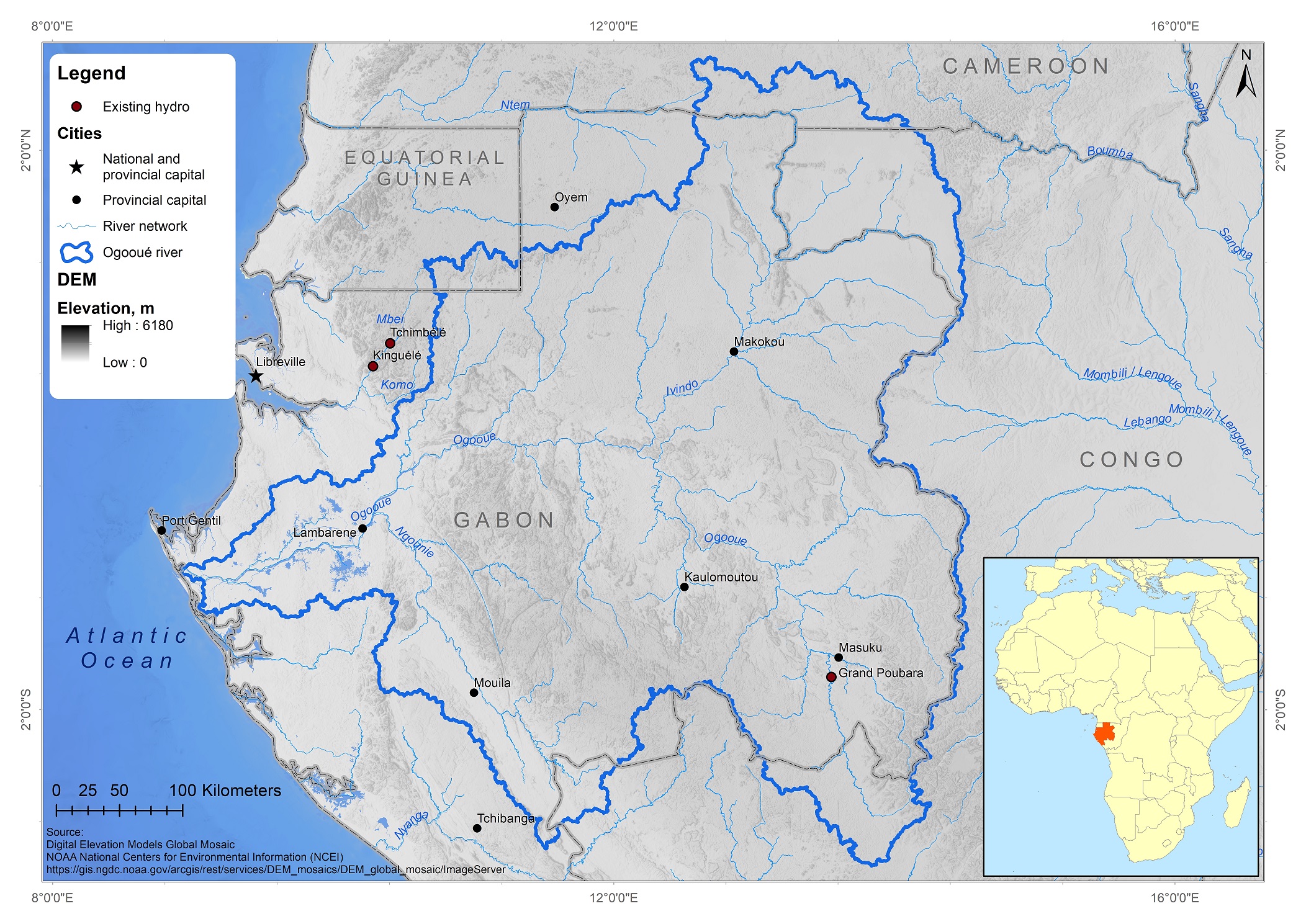}
    \caption{The Ogooué river basin.}
    \label{fig:figure1}
\end{figure}

Gabon’s existing hydropower development consists of three plants in operation. Two of these are in the Mbei river, in the Komo river basin: Tchimbélé (68 MW) and Kinguélé (58 MW), supplying Libreville, the capital and most populated region of Gabon. The third is the Grand Poubara dam (157 MW) located on the Ogooué river. It started operation in 2013 to supply the Franceville system, the third largest urban center in Gabon in the easter part of the country.

Gabon's electrification master plan \cite{27} presents 38 potential hydropower dam sites, with 28 of these located in the Ogooué watershed. Despite being a recent study, hydropower candidates were identified in surveys that date back to the 1980s, when concerns about social and environmental issues caused by the projects were of lesser importance. In this context, updating the hydropower potential survey of Gabon is a critical step to the identification of sustainable hydro candidate projects that can meet this country's growing demand.

Although the Mbei river is very important to the Libreville electrical system, it is already fragmented by two existing dams. The Komo is a free-flowing river, but that would change with the construction of the proposed Ngoulmedjim dam, a scheme that would divert water from the Komo and create a 70 km stretch of reduced flows. Ecological flows have not yet been studied or proposed to the best of our knowledge, which raises a concern about the sustainability of the project. Thus, other than decisions relating to the Ngoulmedjim dam and its ecological flows, there is not a great deal of space for sustainability planning in the Komo basin.

The Nyanga basin, near the border with the Republic of Congo, has a small hydropotential when compared to the other rivers. It is currently undammed and has 679 km of highly suitable fish habitats \cite{28}. Connectivity loss would be large if the Igotchi dam with 28 MW of installed capacity identified in \cite{27} is developed. In this case 392 km of highly suitable habitats would be lost, or 58\% of the total length \cite{28}. This is an example of a project with significant impacts to fish habitants that has a small contribution for increasing the supply of electricity and should probably not be part of a set of projects selected according to sustainable principles.

Thus, it is apparent that the Ogooue river basin is the critical focal point for sustainability planning. It has the largest potential for development, with numerous possible sites. And it has critically important socio-environmental values in its present largely free-flowing state. It is here that sustainability planning has the potential for making the greatest difference in terms of energy production and social environmental conservation. And it is here that we can best illustrate the use of the HERA methodology.

\section{Methodology}
\label{2}

\subsection{Study area description}
\label{2.1}

The Ogooué is the second major river in central Africa, only surpassed by the Congo River in terms of discharge \cite{29}, with a mean flow of 4,700 m³/s. The Ogooué River Basin (ORB) has an area of 224,000 km², of which 90\% in Gabon and the remaining 10\% in Cameroon and Congo-Brazzaville. ORB covers more than 80\% of Gabon’s total area and has several national and international areas of great importance, such as national parks (Lopé National Park, Ivindo National Park, etc.), marine-protected areas, UNESCO’s Natural Heritage site (Portes d’Okanda) and RAMSAR sites (the Rapids and Chutes of the Ivindo, or the Mboungou Baduma and the Doumé Rapids) \cite{26, 27, 28, 29, 30}. The ORB is one of the most preserved ecosystems in Central Africa, almost completely covered with dense vegetation of the central African rainforest. 

ORB is characterized by plateaus and hills bordering a narrow coastal plain \cite{29}. There are three main topographical structures in the basin \cite{31}: (1) the low coastal plain extends from the Ogooué mouth to the Atlantic Ocean and covers the lakes region and the river delta; (2) the plateaus cover most of the surface area of the basin, extend northwards over Ogooué-Ivindo, southeast over Haut-Ogooué (the Batéké plateau), and southwards over the Ogooué-Lolo and Ngounié regions; and (3) the mountainous massif inside the basin constituted by medium-sized mountains.

Although the hills are not very high (mean elevation in the catchment is 450 m), the Ogooué is unnavigable between Lastoursville and Ndjolé due to chutes and rapids. After Ndjolé, the river runs west and reaches the 100 km wide and 100 km long Ogooué delta. The lower part of the Ogooué is navigable and gentler than the rest of the river, with relatively low bed slopes, between 0.07 and 0.13 m km-1 \cite{29}. Its largest tributaries are the Ivindo and the Ngounié.

\subsection{Data sources}
\label{2.2}

The hydrologically adjusted elevations from MERIT Hydro \cite{32} were used in this study. It is a global flow direction map at 3 arc-second resolution (~90 m at the equator) derived from MERIT Digital Elevation Model (DEM) \cite{33} and water body datasets, such as G1WBM, OpenStreetMap, and GSWO. In MERIT Hydro, elevations are also adjusted to satisfy the condition “downstream elevation cannot be higher than upstream” while minimizing the required modifications from the original MERIT DEM. The MERIT DEM, in turn, was developed by removing multiple error components (absolute bias, stripe noise, speckle noise, and tree height bias) from existing spaceborne DEMs, such SRTM3 v2.1 and AW3D-30m v1.

ORB climate is equatorial, with two rain seasons: February to May and October to December. Mean annual precipitation is 1831mm and temperatures vary between 21 and 28 °C \cite{31}. Daily river discharges for gauging stations in Gabon come from SIEREM \cite{34}. They are available for several locations in the basin, and their records in general vary between the 1950s and the 1980s.

The 7Q10 – the lowest 7-day average flow that occurs on average once every 10 years – was used as ecological flow to be maintained in the natural stream between the water intake at the dam and the powerhouse in diversion schemes. This criterion should be reviewed at later stages of the development, after investigations of adequate ecological flows in each location based on local conditions \cite{35}.

For the analysis, we have generated hydropower development alternatives, each of which is a set of selected projects from a candidate list in a river basin, in this case the ORB. The HERA methodology compares these alternatives using impact metrics, designed to measure the interference of reservoirs with spatial layers of information. The sources of the spatial layers and the related impact metrics used in this study are presented in ~\ref{A}.

Some metrics can be cumulative for each alternative. For example, the reservoir flooded area of an alternative is the sum of the flooded areas of the reservoirs in that alternative. On the other hand, river connectivity (or its contrary, fragmentation) is a metric valid for the alternative as a whole and not for the sum of all reservoirs. For example, a dam located lower in the basin would cause fragmentation throughout all upstream segments in the basin, while an additional upstream dam would not cause any additional fragmentation.

\begin{figure}[hbt!]
    \centering
    \includegraphics[width=1\textwidth]{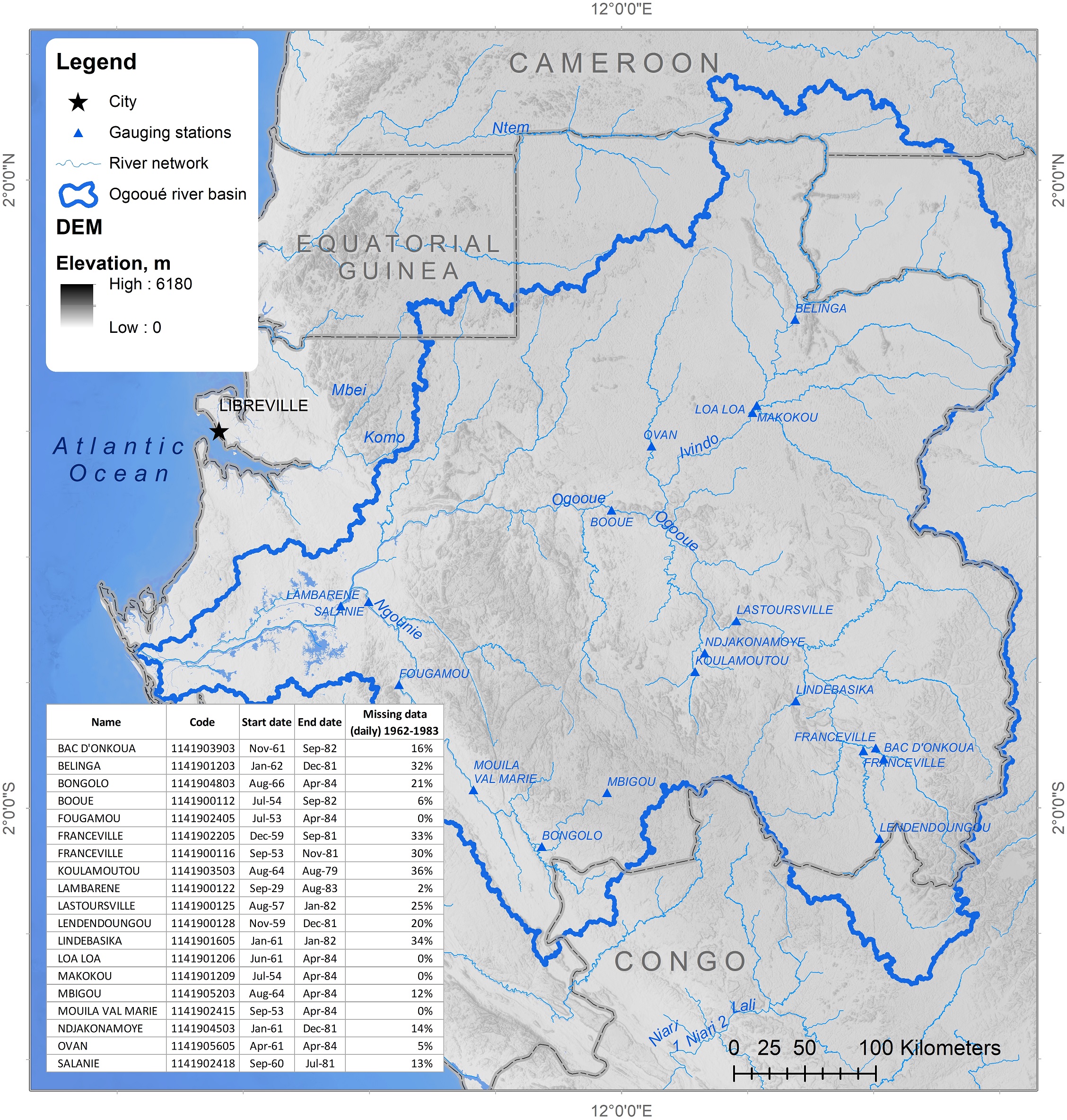}
    \caption{River discharge gauging stations in the Ogooué basin.}
\end{figure}

These metrics can be used as constraints. For instance, the maximum flooded area should be smaller than a provided threshold. They can also be cost components in the economic evaluation of alternatives (e.g. costs related to land acquisition, community resettlement or reconstruction of flooded infrastructure).

One of the challenges to evaluate impacts on freshwater life is the lack of comprehensive studies on fish species and their behavior in Gabon. To overcome this, we rely on studies that used a Maximum Entropy (MaxEnt) distribution model to assess freshwater fish species richness at the landscape level \cite{28,36}. By integrating the findings of these studies with a method for effectively designing new proposed dams in Gabon, we can stipulate project local impacts and those related to the river fragmentation in highly diverse areas.

The length of free-flowing rivers depends on the identified drainage network. In this study the drainage network was identified based on MERIT Hydro DEM, from Ogooué’s mouth to headwaters, excluding the fragmented stretch by the existing Poubara dam.

\begin{figure}[hbt!]
     \centering
     \begin{subfigure}[b]{0.48\textwidth}
         \centering
         \includegraphics[width=\textwidth]{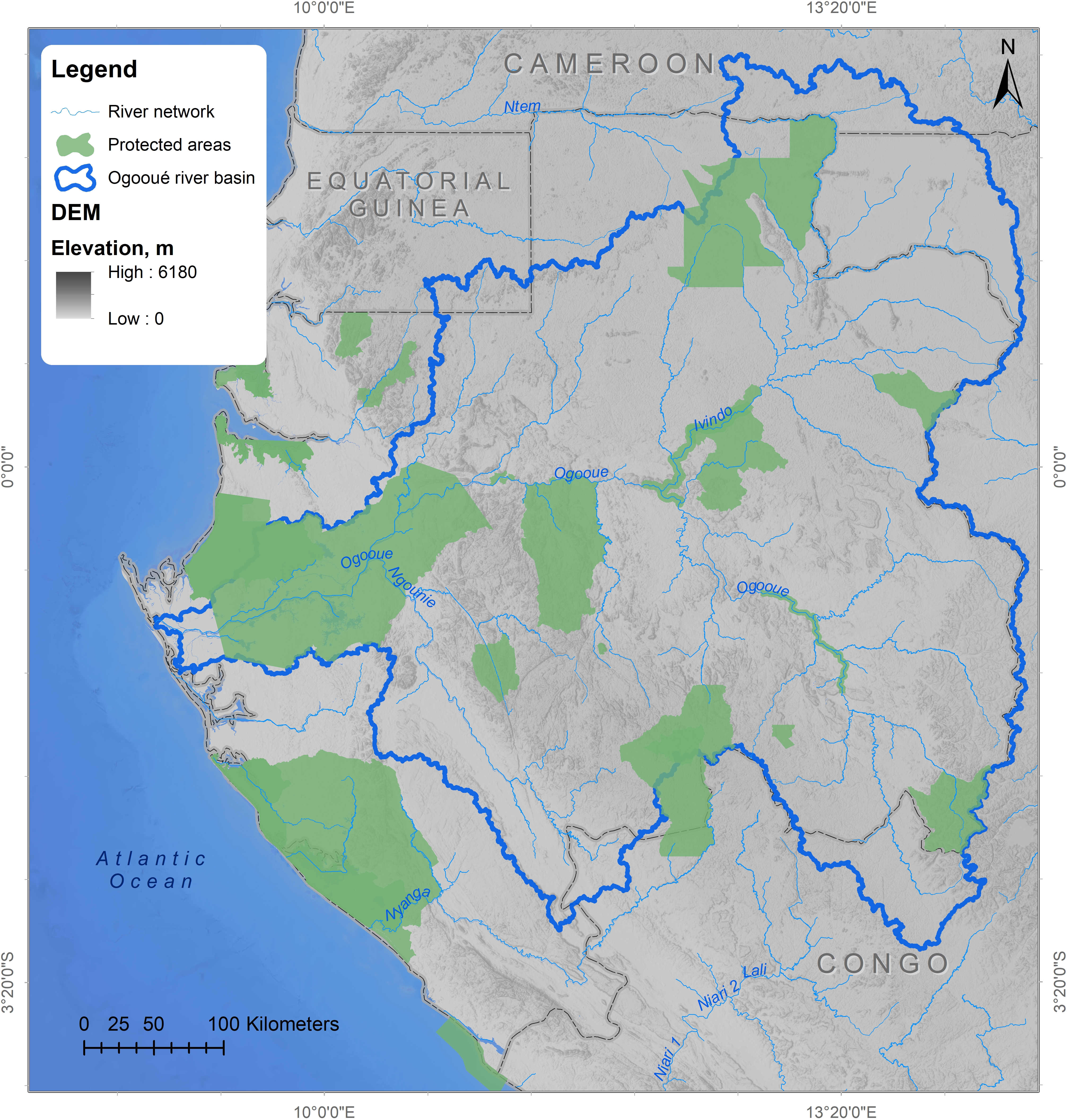}
         \label{Figure 3a}
     \end{subfigure}
     \hfill
     \begin{subfigure}[b]{0.48\textwidth}
         \centering
         \includegraphics[width=\textwidth]{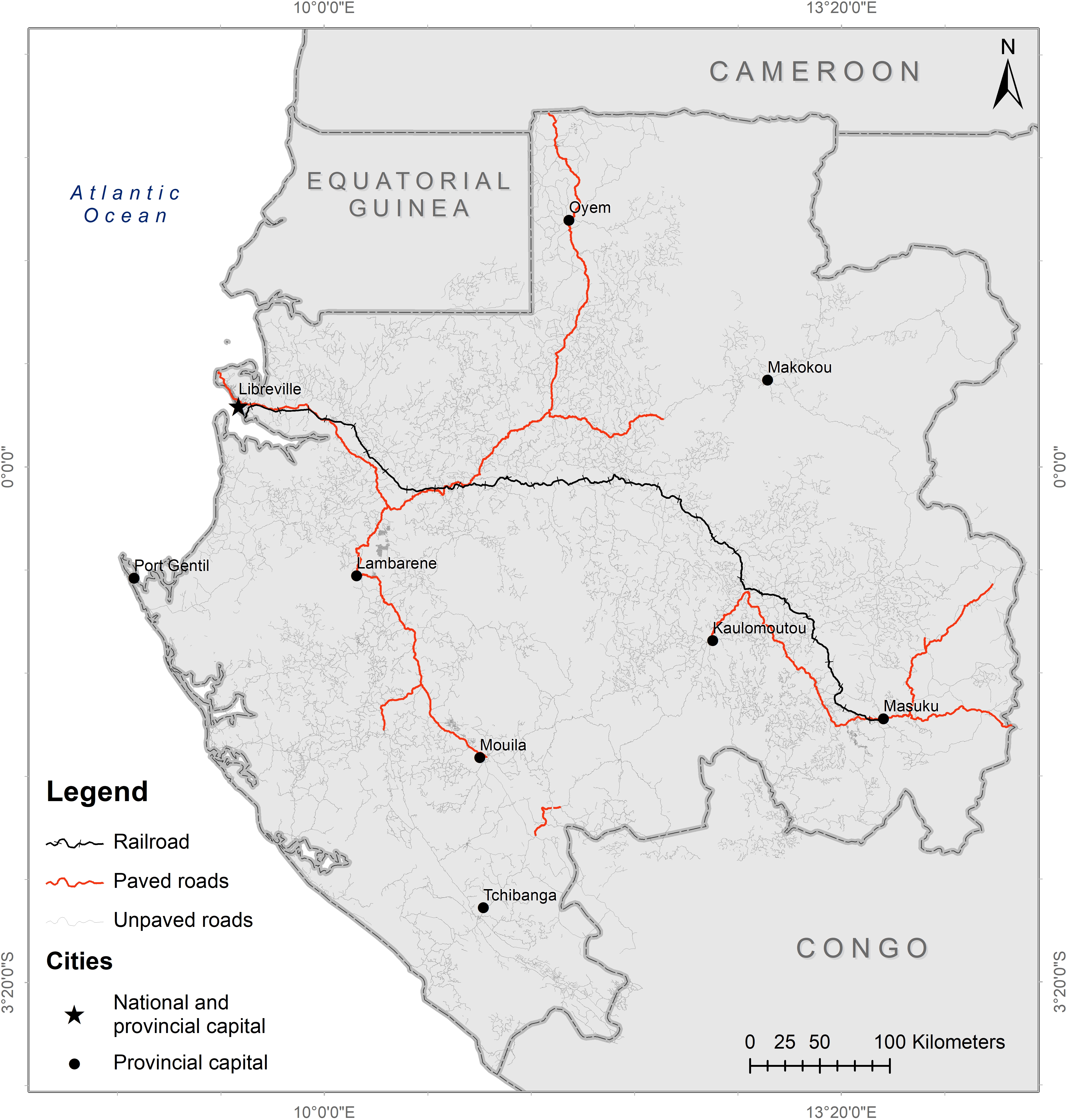}
         \label{Figure3b}
     \end{subfigure}
     \hfill
        \caption{Examples of layers used in the calculation of metrics: protected areas (left) and roads and railroads (right)}
        \label{Figure 3c}
\end{figure}

\begin{figure}[hbt!]
    \centering
    \includegraphics[width=1\textwidth]{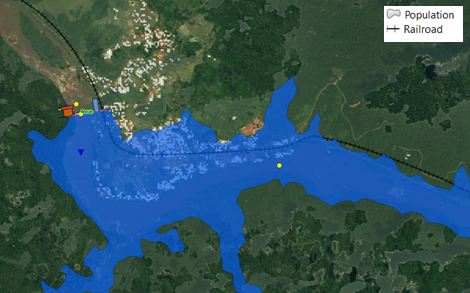}
    \caption{Example of a simulated reservoir and impacts to population and railroads (layers of information).}
    \label{Figure 4}
\end{figure}

\subsection{HERA Components}
\label{2.3}

The framework for sustainable hydropower planning in hydrographic basins applied in this study, named HERA (Hydropower and Environmental Resource Assessment) has three main components, described in the sequence.

\subparagraph{i. Geographic Information System (GIS) processing}\hfill \break
The ORB network is extracted from MERIT Hydro \cite{32}. Historical records of monthly mean and daily maximum river discharges for each year at gauging stations from ORSTOM \cite{34} were used to estimate the discharges at any point of the ORB based on the ratio between the drainage area of the point and corresponding area of the gauging stations. Candidate sites are selected from the ORB drainage network. Each candidate site may have several different hydropower projects, considering different possible water heads. Naturally, at most one project may be selected in each site. This logical constraint is included in the mathematical programming problem formulation. 

\subparagraph{ii. Engineering}\hfill \break
In each site and for each head, an engineering component analyses several engineering arrangements, considering the combination of possible structures (e.g. concrete or earthfill dams, Francis or Kaplan turbine, river diversion or foot-of-dam schemes, ski-jump spillways or stilling basins) and the positioning of hydraulic structures along the dam axis, such as spillway in the right bank and water intake in the left bank or vice-versa). A computational workflow was designed to apply filters and reduce unnecessary processing. For example, a project with 30m of hydraulic head can either use Francis or Kaplan turbines, but not Pelton turbines that require a much larger hydraulic heads (150m or more). An engineering workflow was also developed to design structures and electromechanical equipment; volumes for civil works are obtained from the geometry of the structures and their intersection with the DEM.

For each project, the historical river discharge data is used as an input to compute the corresponding electricity production for a given installed capacity. An iterative method is used to find the optimum installed capacity, that is, the point where the cost of the next MW becomes higher than its economic benefit. 

Budgets are estimated for each candidate with appropriate engineering guidelines. In this case we have used guidelines taken from the Hydropower Inventory Manual of Brazil \cite{37}. For each structure, different volumes (e.g. concrete, reinforcements, rock and earth excavations) are computed and the corresponding costs obtained from a database of unit prices applicable in the region (Gabon).

In addition to civil works and electromechanical equipment, cost components are included as a compensation for socio-environmental impacts, such as relocated population. Some kinds of impacts cannot be addressed by compensations, and in these cases, alternatives should be compared in terms of metrics to minimize interferences and related impacts.

With HERA, this process can be repeated for a large menu of hydropower candidate projects. Computational speedup is possible from distributed processing.

\subparagraph{iii. Optimization}\hfill \break
A mathematical programming (MP) problem for the selection of the best alternative of projects of the river basin is set up. Binary variables are used for the selection of candidate projects. See ~\ref{D} for the details. The objective function maximizes the economic benefit of a set of projects (main decision variables) given by the difference between electricity sales and annualized investment costs taken from the engineering component. The problem includes physical constraints, such as water balance constraints, minimum and maximum storage, minimum and maximum flows through the turbine, and logical constraints (if a reservoir is built in a location, neighboring projects that would interfere with this selection cannot be selected). As mentioned, socio-environmental constraints related to previously defined project metrics can also be defined, such as a maximum number of affected households or a maximum length of roads to be relocated after flooding. Hydrologic variability is considered in the problem through inflow scenarios (historical records). 

\subsection{Application in the Ogooué basin}
\label{2.4}

The proposed HERA methodology was implemented using the HERA computational model. Input data is loaded in the model and an automatic process screens the ORB drainage network for locations that match certain criteria, such as minimum slope or capacity, minimum hydraulic head, maximum distance to roads or to the transmission network. 
River stretches that have met the criteria used in the screening process are used to find potential sites for the hydro development. After the automated screening process is completed, potential sites are reviewed manually using auxiliary functions, such as river profiles and contour lines. These functions help to identify suitable locations for the development of projects (waterfalls, rapids, and narrow valleys).
While the selection of best locations for dams is manually made by using the best engineering judgment, HERA provides computer support in terms of automated engineering design and cost calculation for many projects. 

\begin{figure}[hbt!]
    \centering
    \includegraphics[width=1\textwidth]{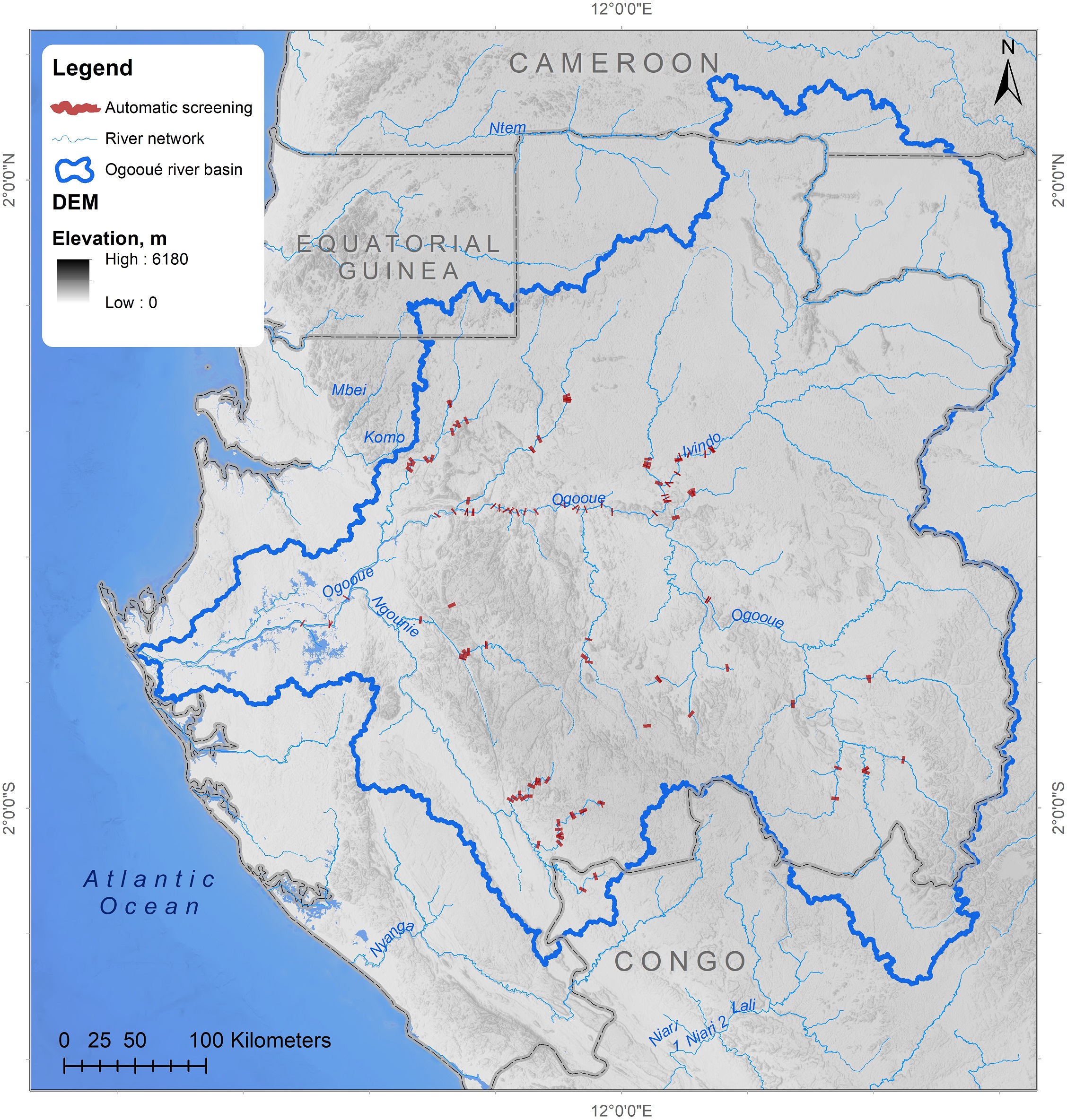}
    \caption{Result of the automatic screening process.}
    \label{Figure 5}
\end{figure}

\begin{figure}[hbt!]
     \centering
     \begin{subfigure}[b]{1\textwidth}
         \centering
         \includegraphics[width=\textwidth]{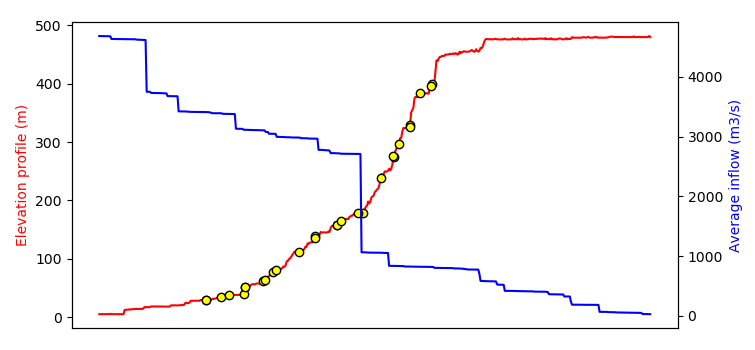}
         \label{Figure 6a}
         \caption{River profile}
     \end{subfigure}
     \hfill
     \begin{subfigure}[b]{0.45\textwidth}
         \centering
         \includegraphics[width=\textwidth]{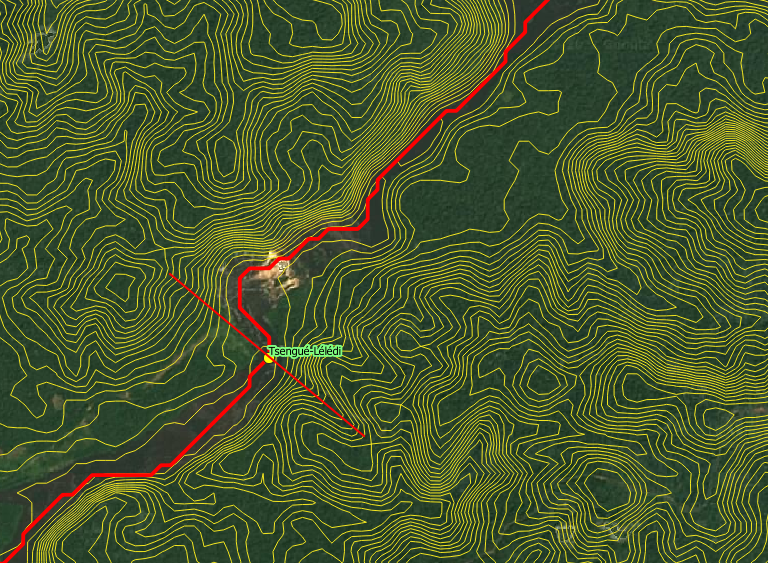}
         \label{Figure 6b}
         \caption{Elevation contour lines}
     \end{subfigure}
     \hfill
       \begin{subfigure}[b]{0.45\textwidth}
         \centering
         \includegraphics[width=\textwidth]{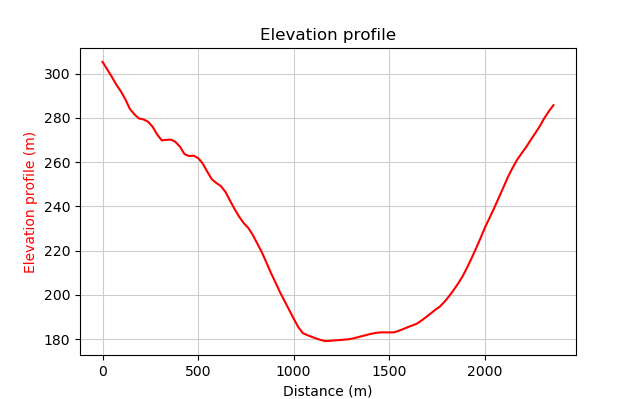}
         \label{Figure 6c}
         \caption{Cross section}
     \end{subfigure}
     \hfill   
        \caption{Example of GIS tools that help identifying potential location for dams and hydro plants.}
        \label{Figure 6d}
\end{figure}

As a result of this screening process, a total of 112 potential sites were identified in the ORB, a large contrast with the 26 projects utilized in the 2017 Master Plan. Considering all hydraulic head and engineering solution variants of each site, a total of 3066 projects were investigated by the engineering component, meaning that their structures were designed, and the final project budget estimated. The list of 3066 projects was filtered according to the following criteria: (i) maximum unit cost 4000 USD per installed kW based upon the assumed power purchase agreement - PPA contract price discussed later in this paper; (ii) minimum power density ratio of 4 MW/km² of flooded area \cite{38}, including riverbed. After applying these criteria, the list of possible projects was reduced to 1539. It is worth mentioning that the list has candidate projects that are mutually exclusive either because they are in the same site or because the construction of one project makes the construction of another unfeasible due to interferences. Appropriate logical constraints are included in the mathematical problem formulation to avoid an infeasible selection of “conflicting” projects.

\begin{figure}[hbt!]
    \centering
    \includegraphics[width=1\textwidth]{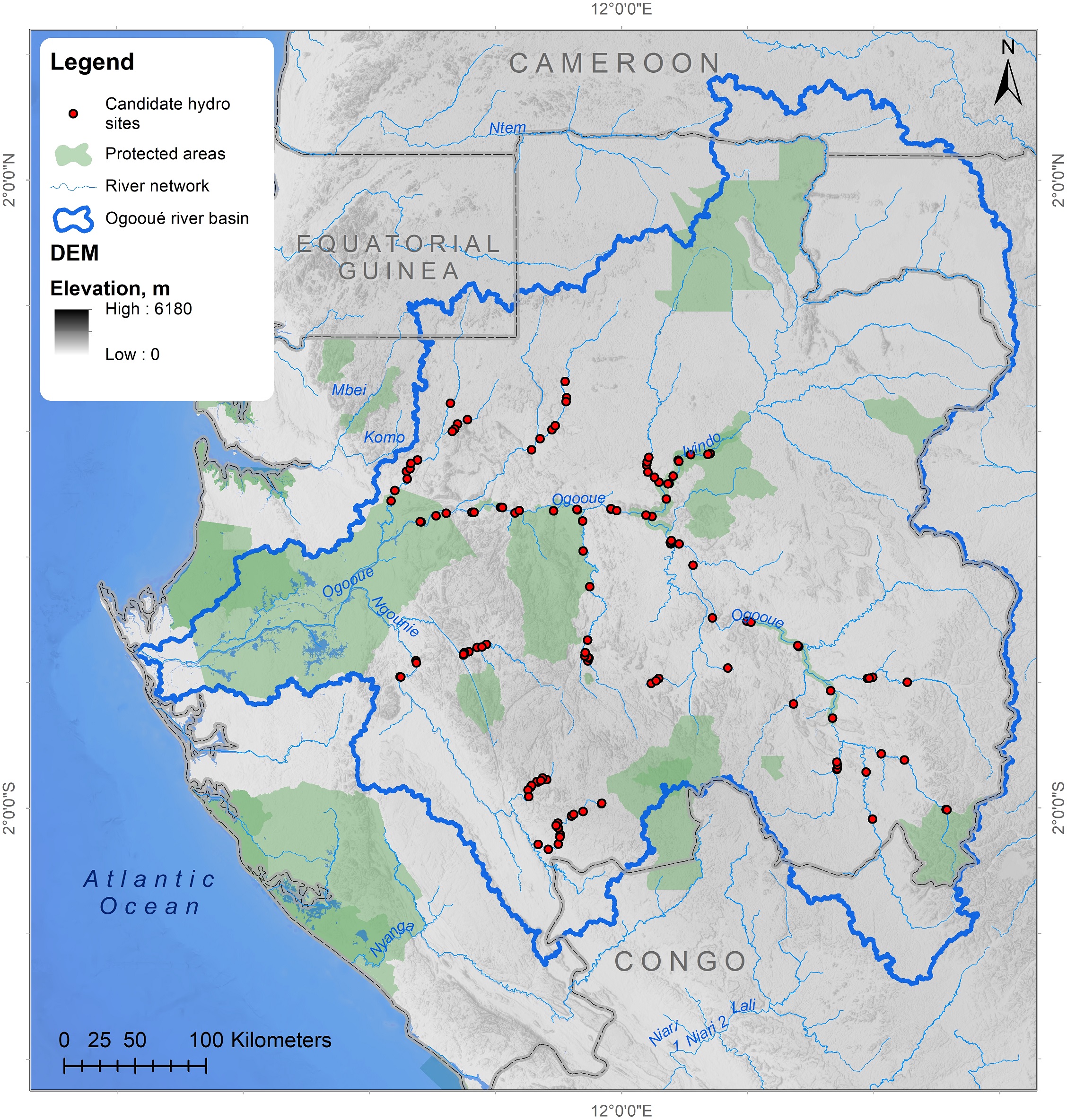}
    \caption{112 identified sites in the Ogooué Basin.}
    \label{Figure 7}
\end{figure}

The following environmental constraints were also considered: (i) the number of relocated households must be less than 100 considering the selection of projects; (ii) reservoirs cannot flood railways; and (iii) no construction allowed in environmentally protected areas.
The yearly revenues of selected hydro projects are their energy production multiplied by an assumed price (USD/MWh) of a Power Purchase Agreement contract. Selected historical inflow years were used with corresponding probabilities. Thus, a sample of yearly energy production or, equivalent PPA revenues is addressed in the problem objective function with total project annuities computed from each project’s CAPEX determined by the engineering module of HERA, a user-provided discount rate (10\% in real terms) and the project useful life (ex. 40 years).

\begin{figure}[!ht]
    \centering
    \includegraphics[width=0.7\textwidth]{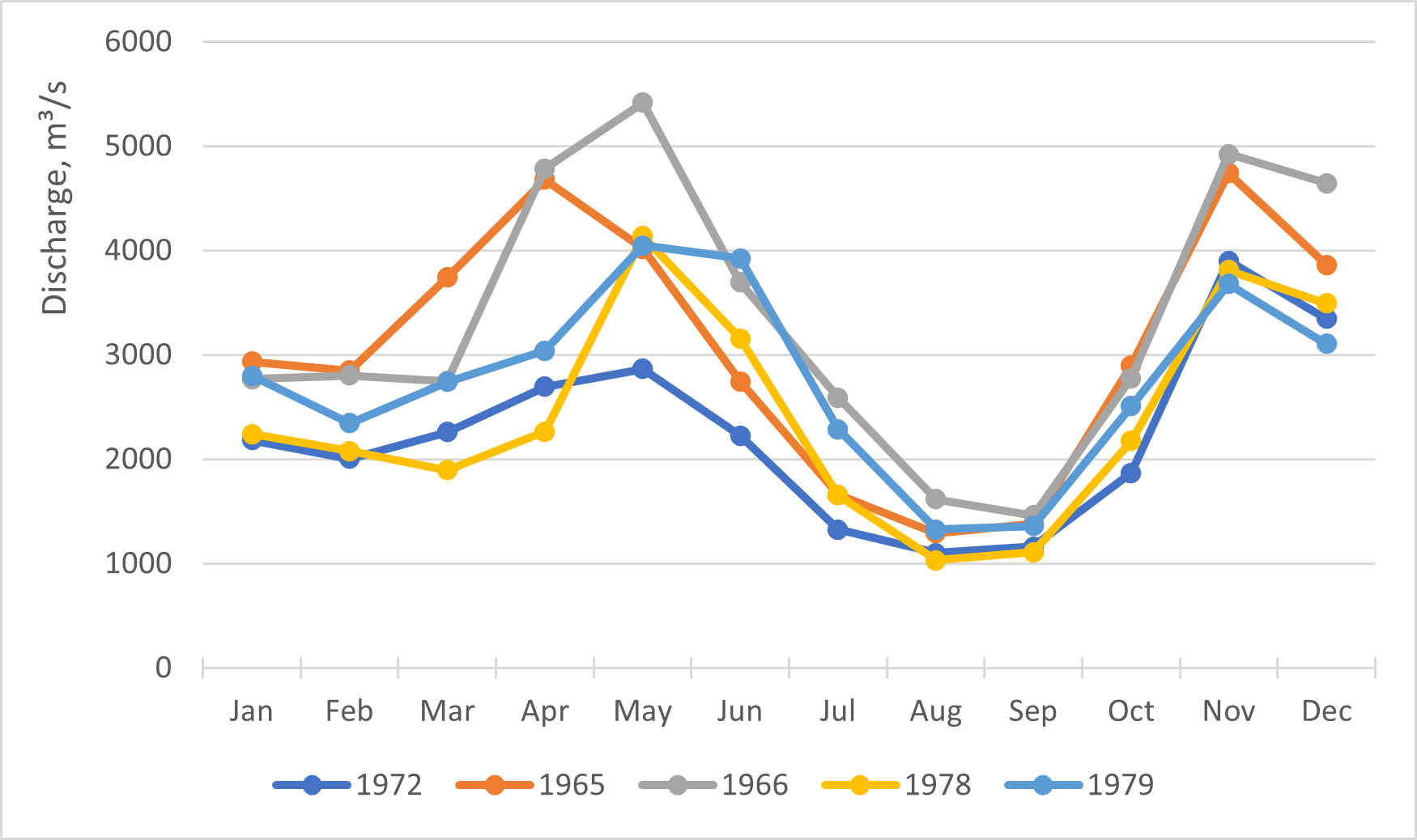}
    \resizebox{3.5cm}{!}{
     \begin{tabular}[b]{cc}\hline
      Year & Probability (\%) \\ \hline
      1972 & 7.4 \\
      1978 & 26.0 \\
      1979 & 33.3 \\
      1965 & 22.2 \\
      1966 & 11.1 \\
      Sum & 100.0 \\ \hline
     \end{tabular}
    }
    \caption{Monthly inflows at Booue for selected hydrological years.}
  \end{figure}

Although hydro plants provide several services to the grid, some related to electricity production, others to the available capacity or ancillary services, in this application we have considered revenues related to PPA sales only.

\section{Results}

After preparing the problem formulation and supplying it with data and parameters from the ORB study case, the resulting mixed-integer (MIP) problem was solved. The optimum alternative included 27 projects that could add 3927 MW to Gabon’s power system, about 5x the total current capacity and well more than the long-term demand of the country. With this selection, however, the length of free-flowing rivers in ORB would be greatly reduced from 15025 km to 4227 km.

Given the abundance of feasible projects we created a risk-adjusted alternative by reducing the price of the energy sold in the PPA. The objective was to select a short list with the most economical projects. We reran the optimization model with a conservatively low PPA price. As a result, only 10 projects that could add 1762 MW to Gabon were selected by the model. The length of free-flowing rivers was reduced to 9150 km, still a significant impact for the country.

Finally, a third alternative was investigated, with the same objective as before, that is, maximize revenues considering a lower PPA price subject to a constraint that imposes a minimum length of free-flowing rivers of 12000 km. The optimization process was finished after finding 20 feasible alternatives. Considering the solution method used \cite{39}, a new feasible solution is only presented if it improves the objective function, thus, total net revenue. Thus, the last solution is inherently the most economically attractive. Nevertheless, it is worth examining in \label{figure9}. the evolution of the optimization process considering not only yearly net revenues (the objective function), but also the number of projects in each feasible solution and other key metrics, such as the aggregate values of installed capacity, free-flowing river length and reservoir flooded area.

\begin{figure}[hbt!]
    \centering
    \includegraphics[width=1\textwidth]{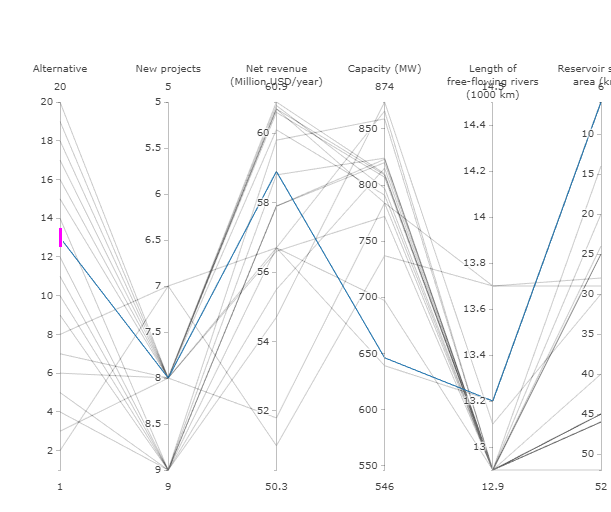}
    \caption{Example of metrics for all possible alternatives generated.}
    \label{fig:Figure 9}
\end{figure}

We notice from the chart that the optimum solution provides 21\% more net revenue than the first solution and only 5\% higher than the 10th solution. By design, all solutions respect the minimum free flowing river length of 12,000km. It can be seen that alternative 13 provides 97\% of the net revenues of the optimum solution, has 13200 km of free-flowing river, but has a much smaller flooded area (only 6 km2 against 46 km2 of the optimum). Considering these additional dimensions, alternative 13 is a well-balanced solution and for this reason is chosen as a representative of this case. It includes 8 projects that can add 646 MW to the power system, roughly the current capacity of the entire country, including thermal power stations.

Alternative 13 presents an interesting additional element. As a result of the additional constraint in fragmentation, the Ivindo sub-basin is left intact. It may be that the intactness of sub-basins is an ecological value worth considering. Given the apparent high species richness \cite{36} in this river, this could be a good conservation alternative.

\begin{figure}[hbt!]
\begin{subfigure}[b]{0.49\textwidth} 
    \includegraphics[width=\linewidth]{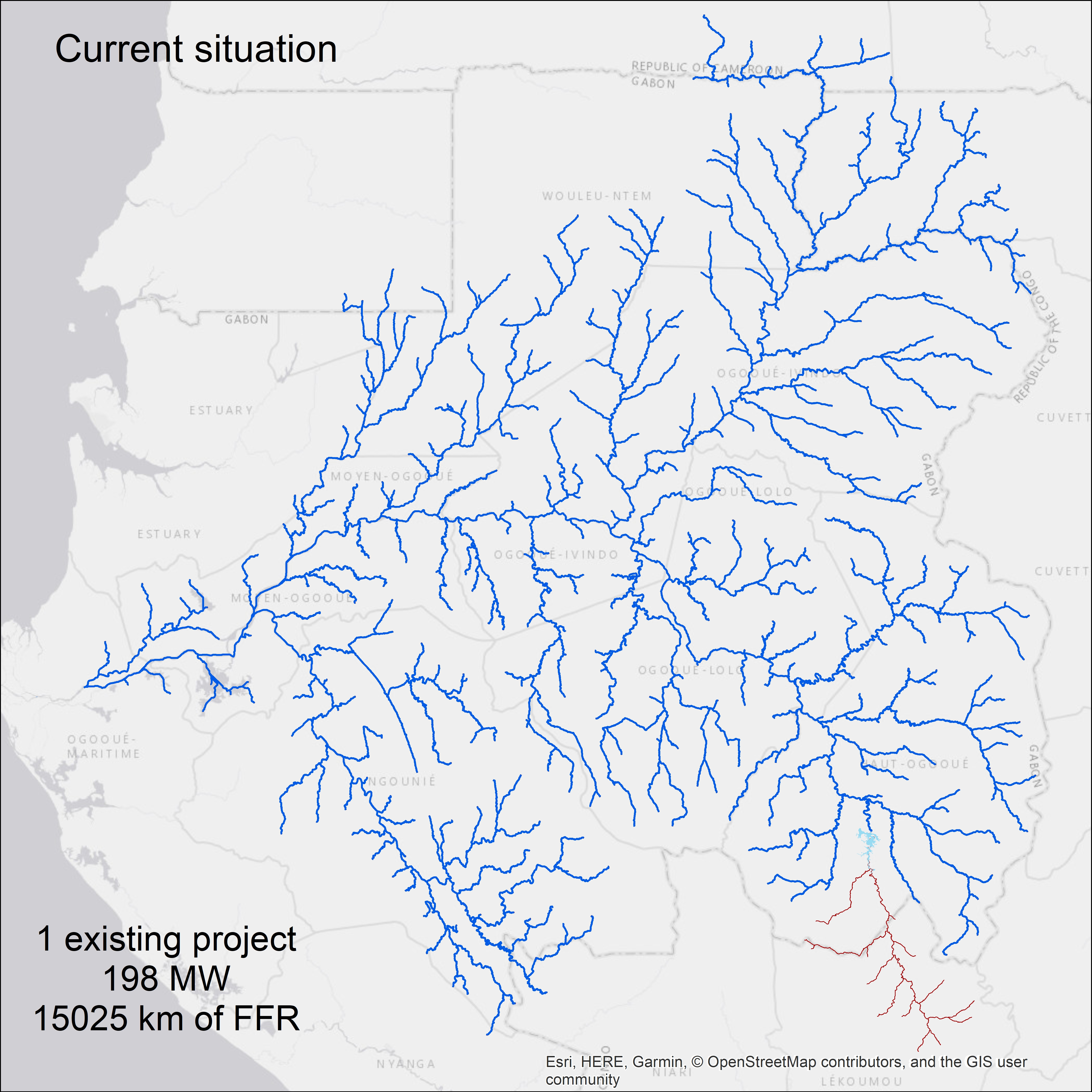}
    \caption{Current situation.} 
\end{subfigure}%
\hfill
\begin{subfigure}[b]{0.49\textwidth} 
    \centering
    \includegraphics[width=\linewidth]{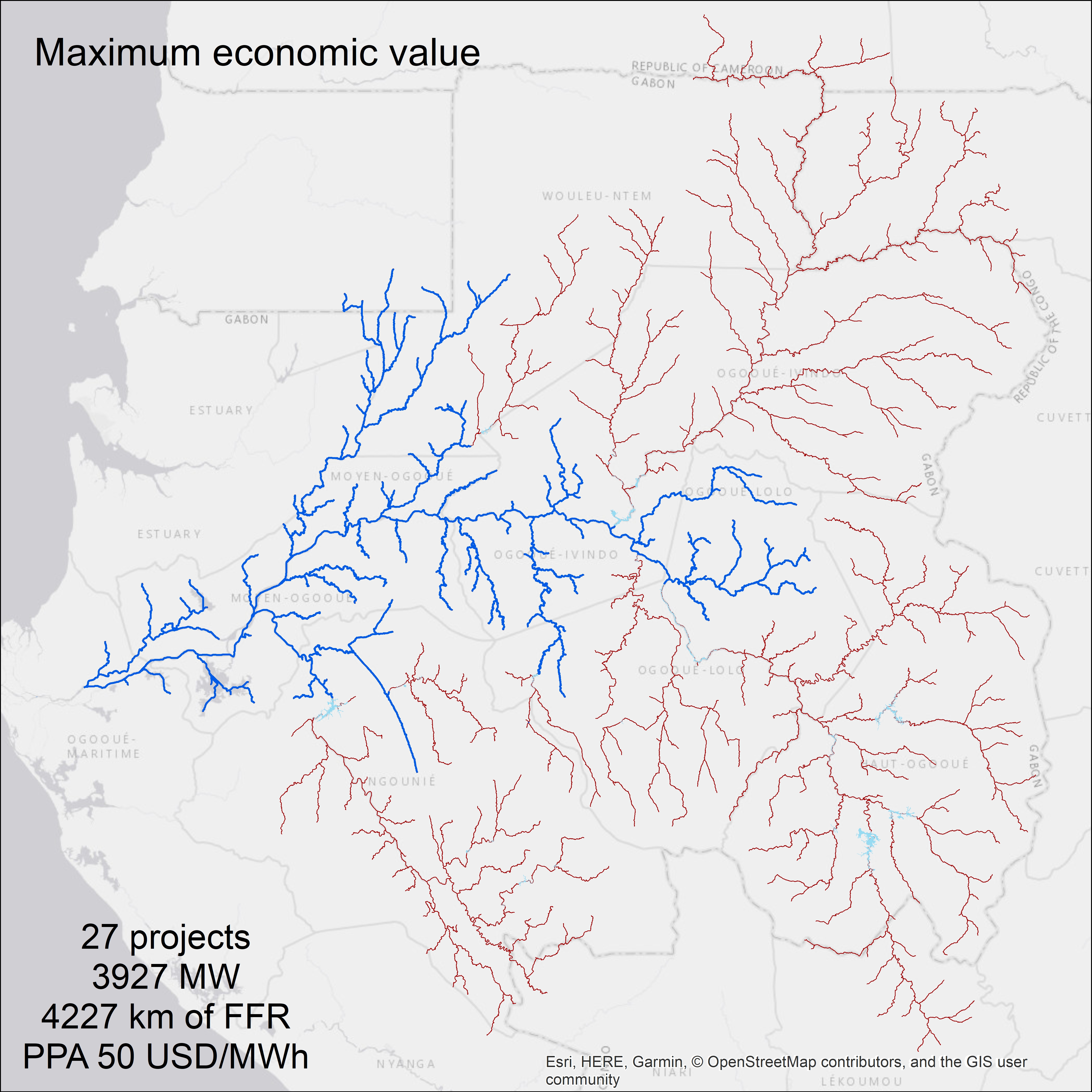}
    \caption{Maximum economic value.} 
\end{subfigure}

\begin{subfigure}[b]{0.49\textwidth} 
    \includegraphics[width=\linewidth]{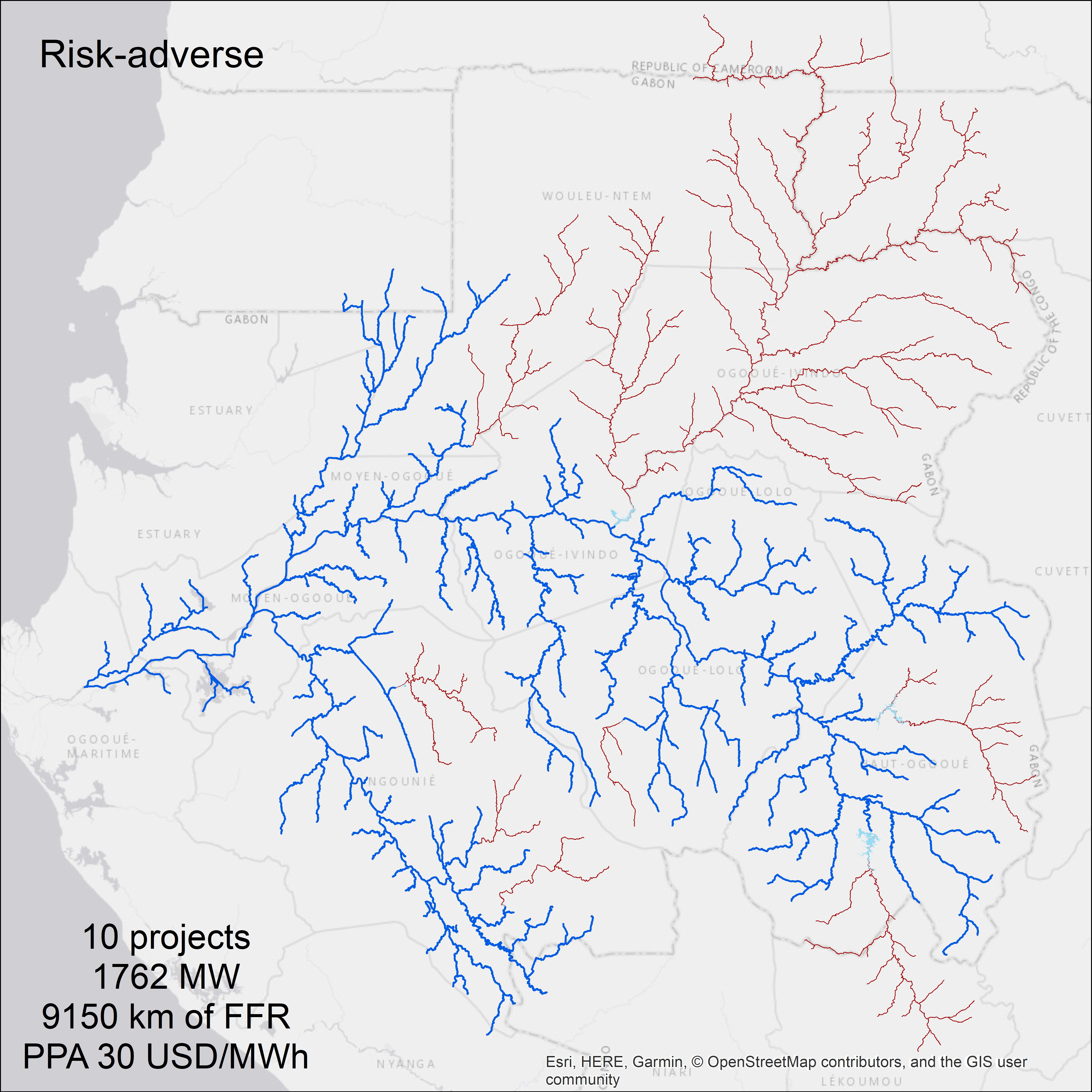}
    \caption{Risk-adverse.} 
\end{subfigure}%
\hfill
\begin{subfigure}[b]{0.49\textwidth} 
    \includegraphics[width=\linewidth]{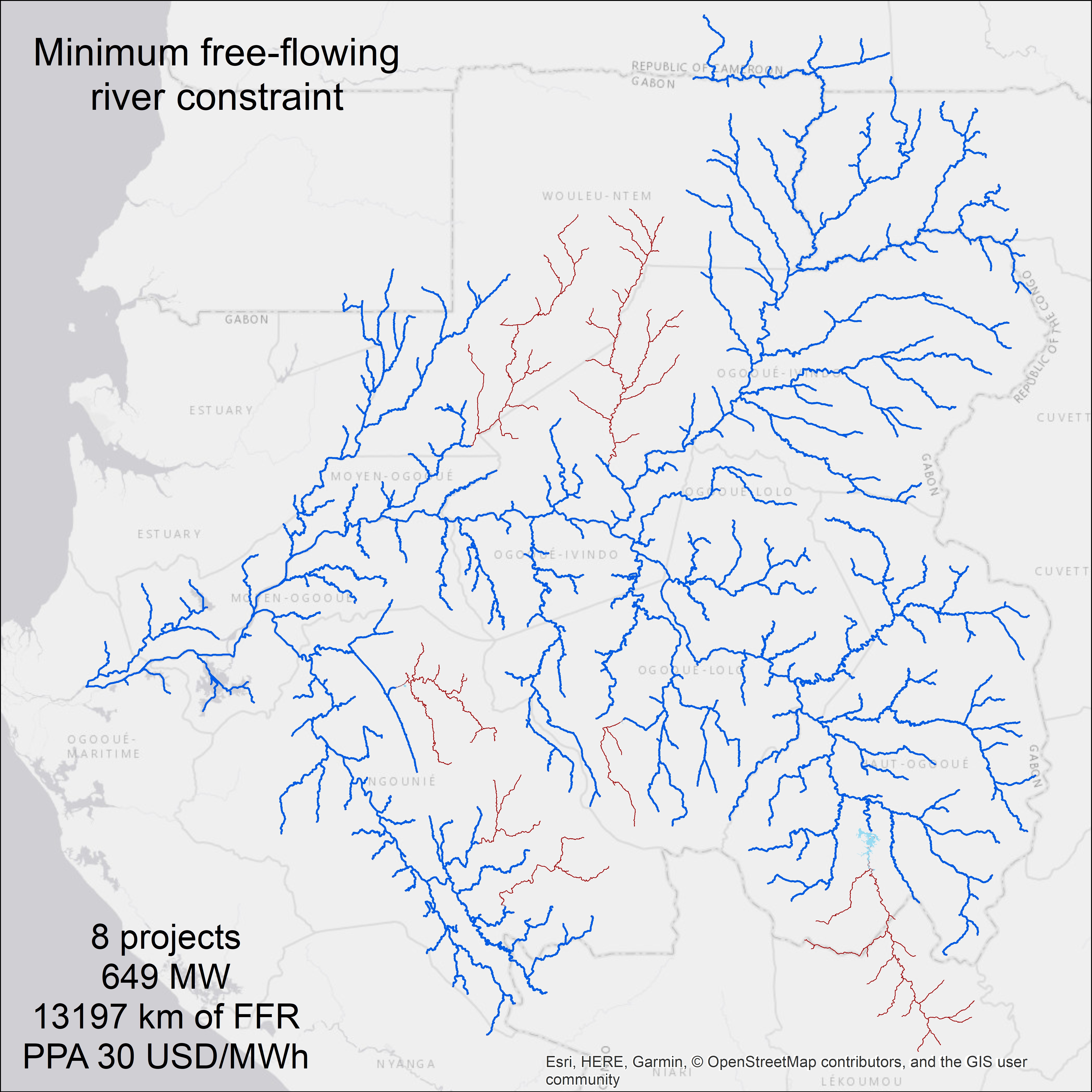}
    \caption{Min. free-flowing river constraint.} 
\end{subfigure}%
    \caption{Fragmentation for the current situation and the 3 analyzed alternatives.} 
\label{Figure10}
\end{figure}

\section{Discussion and conclusions}
\label{3}

We have presented a framework to guide sustainable hydropower development based on the maximization of the economic value of an alternative while limiting impacts through user-defined metrics and optimization constraints. In addition to the optimum alternative, feasible alternatives provide with these metrics an interesting and objective opportunity for stakeholders to discuss possible hydropower development plans, including metrics that were not part of the problem formulation (e.g. reservoir flooded area in the example; or intactness of integral sub-basins). 

A key contribution of the proposed approach is the preparation of thousands of candidate projects and the assessment of their impacts. To achieve success in sustainable hydropower development it is critical to bridge together a multi-disciplinary team that will a holistic view of pros and cons of damming a river. More broadly the exercise with the ORB provided insights to a possible stakeholder engagement covering a three-step process, with specific objectives, as listed below:

\begin{enumerate} 
 \item HERA methodology: 
  \begin{enumerate}
   \item Identify promising locations with screening tools.
   \item Simulate the construction of hydropower projects with proven guidelines for engineering design and estimate corresponding budget.
   \item Apply basic filters ex-ante prior to optimization to remove noncompetitive projects or those with large negative impacts.
   \item Formulate (update formulation of) a model and solve it considering the list of candidate projects selected in the previous step.
  \end{enumerate}

 \item Optimization: 
  \begin{enumerate}
   \item Solve a problem for the list of candidate projects.
   \item Recover sub-optimal feasible solutions (set of selected projects with net revenues sufficiently close to the optimum solution).
  \end{enumerate}

 \item Discussion and validation:
  \begin{enumerate}
   \item Calculate relevant metrics for retrieved solutions.
   \item Evaluate feasible alternatives considering multiple dimensions (metrics) and problem solution (net revenues). In this case, if a consensus is built over design plan, done. If not, evaluate from the incumbent solution metrics whose values exceed a threshold agreed by the stakeholders.
   \item Add constraints to the problem formulation, such as the minimum free-flowing river length used in the case of ORB. Go to step 2 and keep track of revenue-loss incurred by constraints when solving the updated problems. Oftentimes, socioenvironmental compensations that can relax or remove some of these constraints can be a more effective economic solution, leading to better development alternatives.
  \end{enumerate}
\end{enumerate}

Steps 1-3 above can hopefully encourage a transparent and participatory process of hydropower development alternatives. It should be preferably used in the early stages of the planning process to guide project siting and design. An added benefit should be an improved communication between stakeholders and local communities. 

\begin{figure}[hbt!]
    \centering
    \includegraphics[width=1\textwidth]{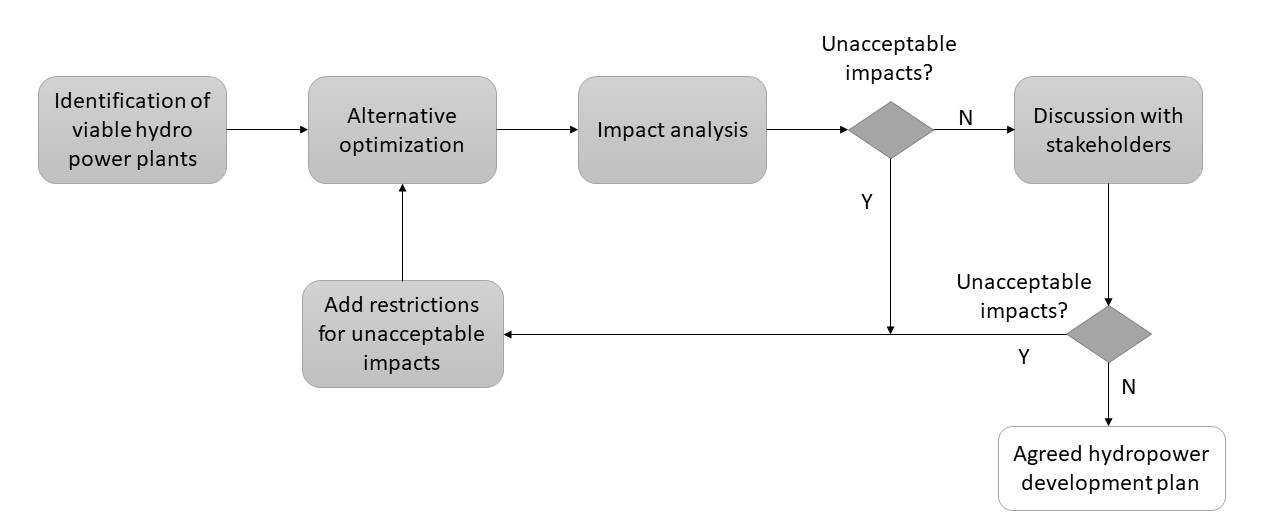}
    \caption{Iterative planning process that can be applied with HERA. }
     \label{Figure 11}
\end{figure}

The interaction between this planning methodology and the integration of variable renewable energy (VRE) sources is a topic for future research considering that the growth of VRE is key for the achievement of net-zero emissions plans. Hydropower can work as a natural battery, supporting electrical grids to expand with renewable supply in a stable and reliable way, without the need to fall back on fossil fuels to avoid blackouts. However, the impact on flow sediment regimes and river discharge alteration must also be considered. Thus, the preparation of hydropower projects (both conventional and pumped-hydro-storage) should be integrated in new studies considering integrated resource planning of power systems. 

\section*{Acknowledgements}
The authors would like to thank David Harrison for his contribution in reviewing this paper. We would also like to thank Marie-Claire Paiz, Erik Martins, Elvis M. Makady and Anne M. Trainor from The Nature Conservancy (TNC) for various contributions that lead to the design of the proposed framework. Finally, we thank the TNC team for the opportunity of working together in hydropower planning Gabon in close collaboration with local governmental institutions.

\appendix

\section{Data sources}
\label{A}

\begin{table}[H]
\centering
\caption{Spatial layers and their sources.}
\begin{adjustbox}{width=\textwidth}
 \begin{tabular}[b]{p{0.3\linewidth} | p{0.7\linewidth}}\hline
      Data & Source \\ \hline
      Aboveground biomass density (Mg/ha) & \url{LC_Nasa_carbon database (https://carbon.jpl.nasa.gov/)} \\
      \\
      Land cover & ESA. Land Cover CCI Product User Guide Version 2. Tech. Rep. (2017). Available at:  \url{maps.elie.ucl.ac.be/CCI/viewer/download/ESACCI-LC-Ph2-PUGv2_2.0.pdf} \\
      \\
      Population & Facebook data \url{https://dataforgood.fb.com/tools/population-density-maps/} \\
      \\
      Parks and protected areas & \url{http://gab.forest-atlas.org/map/?l=en} \\
      \\
      Intact Forest Landscape (IFL) & \url{http://intactforests.org/data.ifl.html} \\
      \\
      Bridges & \url{https://data.amerigeoss.org/id/dataset/02603144-e4af-44f8-86b5-b9f367d40a8e} \\
      \\
      Fish barriers, Mining areas, Forest concessions, Transmission lines, Agriculture, Roads and railroads & Atlas d'Eaux Douces du Gabon \url{https://hub.arcgis.com/maps/e0aeb96fd4304016a1c54f51e647fd51/about} \\ \hline
 \end{tabular}
\end{adjustbox}
\end{table}

\section{HERA templates and layouts}
\label{B}

\begin{figure}[H]
    \centering
    \includegraphics[width=1\textwidth]{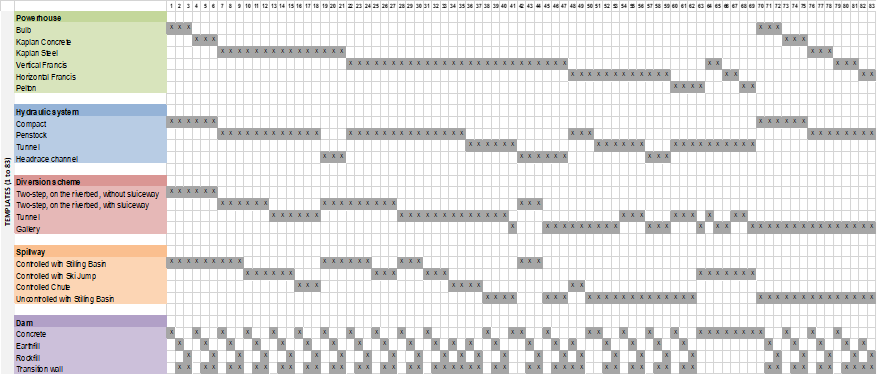}
    \caption{Templates: set of alternatives with 83 combinations of structures. Each column corresponds to a different template. The structures that make part of each template are highlighted in gray with an “x”.}
\end{figure}

\begin{figure}[H]
    \centering
    \includegraphics[width=1\textwidth]{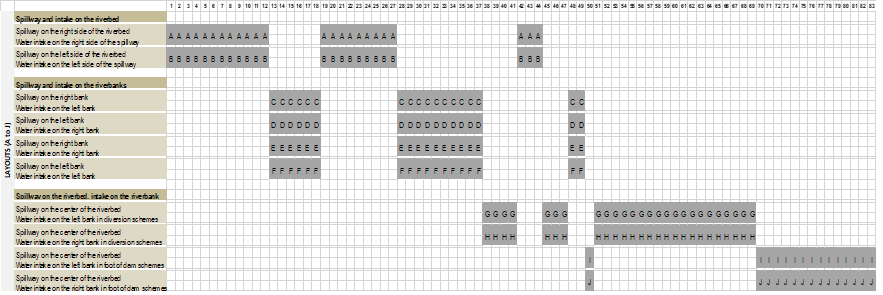}
    \caption{Layouts are alternatives for the positioning of the structures along the dam axis. In each line choose only one column that is filled with a letter, aiming to form a sequence of type A B C D E F G H I.}
\end{figure} \break

\section{Unit prices}
\label{C}
\begin{table}[H]
\centering
\caption{Unit prices for civil works.}
\begin{adjustbox}{width=\textwidth}
\begin{tabular}[b]{p{0.75\linewidth} | p{0.05\linewidth} | p{0.20\linewidth}}\hline
Item &	Unit &	Cost (US\$ per Unit) \\ \hline
Common excavation &	m³ &	9 \\
Surface rock excavation &	m³ &	24 \\
Underground rock excavation &	m³ &	139 \\
Borrow soil &	m³ &	12 \\
Quarry rock &	m³ &	34 \\
Foundation cleaning and treatment - dam earthworks &	m² &	8 \\
Foundation cleaning and treatment - concrete structures &	m² &	37 \\
Cofferdam removal &	m³ &	10 \\
Cofferdam - 1st phase &	m³ &	9 \\
Cofferdam - 2nd phase &	m³ &	9 \\
Compacted earth fill &	m³ &	8 \\
Clay core &	m³ &	11 \\
Rockfill &	m³ &	11 \\
Filters and transitions &	m³ &	31 \\
Riprap or rockfill protection &	m³ &	24 \\
Downstream face protection (grass) &	m² &	15 \\
Cement &	ton &	180 \\
Structural concrete &	m³ &	220 \\
Roller-compacted or mass concrete &	m³ &	120 \\
Shotcrete &	m³ &	300 \\
Reinforcement steel &	ton &	2,800 \\
Steel lining of penstocks &	ton &	8,000 \\
\end{tabular}
\end{adjustbox}
\end{table}

\begin{table}[H]
\centering
\caption{Unit prices for environmental and social related items.}
\begin{adjustbox}{width=\textwidth}
\begin{tabular}[b]{p{0.4\linewidth} | p{0.3\linewidth} | p{0.3\linewidth}}\hline
Item &	Unit &	Cost (US\$ per Unit) \\ \hline
Road relocation &	m &	750 \\
Railway relocation &	m &	3,000 \\
Bridge relocation &	m &	71,500 \\
Resettlement in rural areas &	household &	15,000 \\
Resettlement in urban areas &	household &	15,000 \\ \hline
\end{tabular}
\end{adjustbox}
\end{table}

The cost estimate of each candidate project also considers indirect costs for miscellaneous items and other items of the budget. They are also informed  by percentages. The following table lists the criteria applied in these cases.

\begin{table}[H]
\centering
\caption{Percentages applied to other items.}
\begin{adjustbox}{width=\textwidth}
\begin{tabular}[b]{p{0.9\linewidth} | p{0.1\linewidth}}\hline
Item &	(\%) \\ \hline
Other costs for civil accounts &	2 \\
Other costs for social and environmental accounts &	30 \\
Miscellaneous items for civil accounts &	20 \\
Miscellaneous items for equipment accounts &	20 \\
Miscellaneous items for social and environmental accounts &	20 \\
Miscellaneous items for indirect costs &	10 \\
Indirect costs for construction site and worker's camp &	2 \\
Indirect costs for maintenance and operation of the site and the camp &	2 \\
Indirect costs for basic engineering &	4 \\
Indirect costs for engineering special services &	1 \\
Indirect costs for environmental projects and studies &	1.5 \\
Indirect costs for owner´s administration &	1.5 \\
Other costs for civil accounts &	2 \\ \hline
\end{tabular}
\end{adjustbox}
\end{table} \break

\section{HERA mathematical formulation}
\label{D}
We now present a novel modeling approach for the optimization of hydropower development that considers the river basin as the planning unit. A mathematical programming formulation is proposed to select candidate hydro power plants from a list of candidate projects such that the present value of net revenues (the difference between electricity revenues and development costs, including construction, electromechanical and socio-environmental) is maximized. The stochastic nature of future water flows also adds to the problem complexity. 
Candidate hydropower projects are created, and their development costs estimated by simulating their construction using different engineering design alternatives based on provided topography. Time series of water inflows to the candidate projects are prepared based on their location, existing gauging stations measurements and geoprocessing functions that can be executed in the clouds for speedup. The modeling approach proposed is quite general and may be used for hydropower assessment in river basins worldwide. It may also be adapted to include other economic uses of water in addition to hydropower. The proposed framework will benefit from the increased availability of topographic data due to technological advances, such as nanosatellites.\\
\\
\textbf{Indices}
\begin{description}
\item[i] candidate projects (total of \textit{I} projects)
\item[t] month
\end{description}
\textbf{Parameters}
\begin{description}
\item $c_i$ Project investment annuity cost (\$/yr)
\item $\rho_i$ production factor of project $(kW/m^3/s)$, related to the project hydraulic head
\item $\overline{\textit{v}_i}$ maximum active storage $(m^3)$
\item $\overline{\textit{u}_i}$ maximum flow thru the turbines $(m^3/s)$
\item \textit{$E_i$} Installed capacity of project (kW)
\item \textit{$h_i$} elevation of the ground at the foot of the dam (meters above sea level)
\item $\delta_i$ hydraulic gross head (m)
\end{description}

Project construction costs are calculated by a specific module which has been integrated to the proposed framework. This module has a built-in “construction engineering logic” and a set of equations for the pre-dimensioning of the project structures (channels, spillways, water intakes, powerhouses, river diversion structures, and others). The result of the module is a sound construction cost estimate $c_i$, which is then input to the mathematical programming (optimization) model and with the following components (1) Basic structures; (2) Dam, spillway, and other civil work construction costs; (3) Turbines and generators costs; (4) Other electrical equipment; (5) Relocations and environmental programs costs; (6) Highways, railroads, and bridges; (7) Other indirect costs. The annuity of the total cost of the project, $c_i$ (\$/year), is computed as the sum of all terms listed above multiplied by $\gamma=(\alpha-1)/(\alpha^T-1)$, where $\alpha=(1/(1+\gamma))$,$\gamma$ is the annual discount rate and \textit{T} is the number of years the project will operate (lifetime).\\
\\
\textbf{Input data}\\
$a_i^t$  series natural of incremental natural inflows to project \textit{i} in month \textit{t} $(m^3/s)$.
$d_t$ duration of month \textit{t} (hours).\\
\\
\textbf{Variables}\\
$x_i$  binary decision variable denoting projects that are selected $(x_i=1)$. If a project is forbidden, $x_i=0$. If a project is obligatory or existing then $x_i=1$.\\
$e_i$ annual electricity production (kWh).\\
$v_{(i,t)}$	water storage in reservoir in month t $(m^3)$.\\
$u_{(i,t)}$	water outflow thru the turbine $(m^3/s)$.\\
$w_{(i,t)}$	spilled outflow $(m^3/s)$.\\
\\
\textbf{Formulation}\\
The mathematical model maximizes the net benefit of the construction of the hydropower projects, which is given by the benefit of the selected projects minus the associated cost. Thus:

\begin{equation} \label{eq1}
Max \sum_i (B_i-c_i \cdot x_i)
\end{equation}

The optimum values of the installed capacity, the corresponding mean electricity and total costs are pre-processed based on the procedures described in this article. The mean production for a given installed capacity is calculated by simulating the project operation with the historical inflows record and with efficiency factors that account for the losses from the conversion of mechanical energy into electrical energy. An unavailability factor is also applied.

The annual benefit of a candidate project \textit{i} is given by the product of the electricity price $\pi_1 (\$/kWh)$ and the yearly electricity production $e_i (kWh)$ sold in a Power Purchase Agreement (PPA) plus the price for capacity payments $\pi_2 (\$/kW/year)$ multiplied by the available capacity $E_i (kW)$. In energy-only markets (no capacity payments), $\pi_2$ is equal to zero.

\begin{equation} \label{eq2}
B_i= \pi_1 \cdot e_i + \pi_2 \cdot E_i \cdot x_i \hspace{50pt} \forall \textit{i} = 1 ... \textit{I}
\end{equation}

The yearly production depends on hydrological conditions, which means there is a natural variability of the PPA revenues. An extension of the model, not presented here for the sake of simplicity, includes inflows scenarios (e.g. different years of inflow measurements) with associated probabilities to estimate the expected production of electricity sold or any measure of energy sold at the PPA, such as the firm energy or any measure of the annual produced energy with given reliability. In the following formulation we are assuming the PPA is remunerating the produced energy of each month and not a statistic of this random variable.

\begin{equation} \label{eq3}
e_i= \sum_{t=1..T} (d_t \cdot \rho_i \cdot u_{i,t}) \hspace{50pt} \forall \textit{i} = 1 ... \textit{I}
\end{equation}

The water balance equation for the reservoir of each selected project states that the initial storage of the next month is the initial storage of the present month \textit{plus} natural incremental volume \textit{minus} the turbined and spilled outflows of the project \textit{plus} those from the set of projects $\Omega_i$ located immediately upstream. The term $3600 \cdot d_t$ is the number of seconds per month \textit{t}. It converts water flow $(m^3/s)$ into storage volumes inflowing to or outflowing from the reservoir in each month.

\begin{equation} \label{eq4}
v_{i,t+1}=v_{i,t} +3600 \cdot d_t \cdot (a_{i,t}-u_{i,t}-w_{i,t}+ \sum_{n\in\Omega_i}(u_{n,t}+w_{n,t}))
\end{equation}

$\hspace{220 pt} \forall \textit{i} = 1 ... \textit{I}, \textit{t} = 1 ... \textit{T}$

A steady state condition is useful to avoid reservoir depletion at end of year $v_{(i,T+1)}$ to maximize turbined flow, thus revenue. We use as boundary conditions that final storage should equate the initial storage, both of which are decision variables.

\begin{equation} \label{eq5}
v_{i,T+1} = v_{i,1} \hspace{50 pt} \forall \textit{i} = 1 ... \textit{I}
\end{equation}

The reservoir storage in each month is limited to the maximum useful storage if the project is built and zero, otherwise. Notice that if $x_i=0$ the right-hand side of the constraint vanishes.

\begin{equation} \label{eq6}
v_{i,t} \leq \overline{v_{i}} \cdot x_{i} \hspace{50pt} \forall \textit{i} = 1 ... \textit{I}, \textit{t} = 1 ... \textit{T}
\end{equation}

The turbines flow in each month is limited to the maximum turbined outflow if the project is built and zero, otherwise. Notice that if $x_i=0$ the right-hand side of the constraint vanishes.

\begin{equation} \label{eq7}
u_{i,t} \leq \overline{u_{i}} \cdot x_{i} \hspace{50pt} \forall \textit{i} = 1 ... \textit{I}, \textit{t} = 1 ... \textit{T}
\end{equation}

Logical constraints are included in the formulation to remove conflicting or mutually exclusive projects, as the next figure shows. Notice that if project \textit{i} is selected, project \textit{k} cannot be selected because the difference between the elevations of project \textit{i} and \textit{k} $(h_k-h_i)$ is smaller than the hydraulic head of project \textit{i}, $\delta_i$. 

\begin{equation} \label{eq8}
x_{i} + x_{k} \leq 1 \hspace{50 pt} \forall \textit{i} = 1 ... \textit{I}, \textit{k} \mid h_k-h_i \leq \delta_i
\end{equation}

Notice that these constraints are also applicable for projects in the same location with different hydraulic heads.

\begin{figure}[H]
    \centering
    \includegraphics[width=1\textwidth]{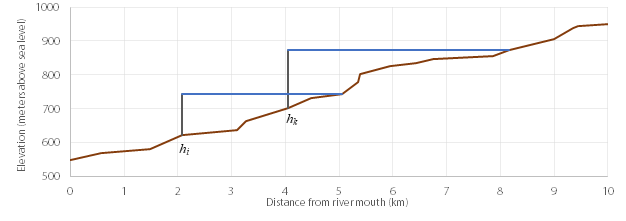}
    \caption{Example of exclusive projects.}
\end{figure} \break

\textbf{Cumulative socioenvironmental constraints}\\
Constraints can also be selected to minimize the cumulative impact of selected projects for given attributes or metrics, such as total flooded area by projects or number of displaced people. Let $m_{i,j}$ be the $j^{th}$ contribution of project \textit{i} to metric \textit{j} with maximum or minimum tolerable values $\overline{M_j}$ or $\underline{M_j}$. These constraints are written as:

\begin{equation} \label{eq9}
\sum_{i} m_ {i,j} \cdot x_i \leq \overline{M_j}
\hspace{20pt} or \hspace{20pt}
\sum_{i} m_ {i,j} \cdot x_i \geq \underline{M_j}
\hspace{50pt}
\forall j = 1...J
\end{equation}

A possible extension is to introduce a \textit{satisfaction} function for each metric \textit{j}, as commonly used in neural networks. Its value is 0, if the metric is less than a minimum acceptable threshold $\underline{M_j}$; and its value is 1, if the attribute is higher than a saturation point $\overline{M_j}$, with intermediate values in between, which is the range of satisfaction improvement of this metric \textit{j}:

\begin{figure}[h]
    \centering
    \includegraphics[width=0.5\textwidth]{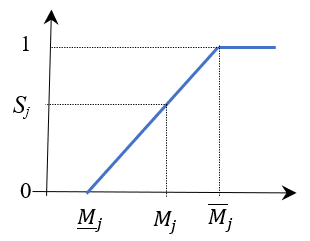}
    \caption{Satisfaction function.}
\end{figure} \break

Mathematically, the satisfaction function is written as:

\begin{equation} \label{eq10}
S_{j} \leq (\sum_i m_ {i,j} \cdot x_i - \underline{M_j}) / (\overline{M_j} - \underline{M_j})
\hspace{50 pt}
\forall j = 1...J
\end{equation}

With these functions HERA’s framework can explicitly include tradeoffs when selecting alternatives. The idea is that an alternative that performs well for a metric and is disastrous for another will be discarded if the bounds for these metrics are included. The proposed approach combines the \textit{mean} satisfaction with the \textit{smallest} satisfaction throughout the metrics, as shown:

\begin{equation} \label{eq11}
\lambda (\sum_{j=1}^{J} S_{j} / J) + (1-\lambda) \cdot \underline{S} \geq S^*
\hspace{50pt}
\forall j = 1...J
\end{equation}\\

Where:\\
$\lambda$	is a weighting factor. If $\lambda =1$, all weight is given to the mean satisfaction; If $\lambda =0$, all weight is given to the \textit{minimum} satisfaction (of all attributes). \\
$S_j$	Satisfaction of metric \textit{j}. \\
$\underline{S}$	the minimum satisfaction of all attributes, that is, $\underline{S} = min{(S_j)} \hspace{5pt} \forall j=1...J$. \\
$S^*$ 	a measure of the combined satisfaction measure.\\

And the decision variable $\underline{S}$ is the minimum value of $min{(S_j)}$ given by \textit{J} inequalities

\begin{equation} \label{eq12}
\underline{S} \leq S_{j}
\hspace{50pt}
\forall j = 1...J
\end{equation}\\

Alternatives produced with this approach have an interesting appeal the application is quite straightforward: stakeholders decide on bounds $\underline{M_j}$ and $\overline{M_j}$ of each attribute, then run the model for various values of $S^*$ and record the output solution (project selection alternative). This approach adds constraints that depend on project selection.\break
\\
\textbf{Minimum river connectivity constraints}

We will simplify notation by assuming that in each candidate sites are there is a single candidate project. We can then write auxiliary binary variables $y_i$  that are activated in river stretches that are fragmented (D.13) by the construction of projects and zero otherwise. “Domino effect” constraints are written in (D.14) such activating auxiliary variables y for all segments upstream of the site where a project was selected. Finally, constraint (D.15) establishes a minimum free-flowing river length $\underline{L}$ which should be met.

\begin{equation} \label{eq13}
y_i \geq x_{i} \hspace{50pt} \forall i=1...I
\end{equation}
\begin{equation} \label{eq14}
y_i \geq y_{D(i)} \hspace{50pt} \forall i=1...I
\end{equation}
\begin{equation} \label{eq15}
\sum_i l_i (1-y_{D(i)}) \geq \underline{L}
\end{equation}

Projects with fish-ladders that avoid river fragmentation or natural barriers can be defined by relaxing constraints (D.13). And free-flowing constraints for subbasins can also be defined. In this later case, set of constraints (D.13) - (D.15) would be replicated in each case with the subset of nodes of the graph of the subbasin.

 \bibliographystyle{elsarticle-num} 
 \bibliography{cas-refs}





\end{document}